\documentclass[12pt,reqno]{amsart}

\usepackage{amssymb}

\renewcommand{\phi}{\varphi}

\newcommand{\Lam}{\Lambda}
\newcommand{\Ome}{\Omega}
\newcommand{\Sig}{\Sigma}

\newcommand{\Z}{{\mathbb Z}}

\newcommand{\seq}{\subseteq}
\newcommand{\stm}{\setminus}

\newcommand{\longc}{,\dotsc,}

\newcommand{\lpr}{\left(}
\newcommand{\rpr}{\right)}

\renewcommand{\>}{\rangle}

\theoremstyle{plain}

\newtheorem{theorem}{Theorem}

\newcommand{\reft}[1]{~\ref{t:#1}}

\newcommand{\refb}[1]{~\cite{b:#1}}

\begin{document}
\baselineskip=16pt

\title{On the size of dissociated bases}

\author{Vsevolod F. Lev}
\address{Department of Mathematics, The University of Haifa at Oranim,
  Tivon 36006, Israel.}
\email{seva@math.haifa.ac.il}

\author{Raphael Yuster}
\address{Department of Mathematics, The University of Haifa,
  Haifa 31905, Israel.}
\email{raphy@math.haifa.ac.il}


\begin{abstract}
We prove that the sizes of the maximal dissociated subsets of a given finite
subset of an abelian group differ by a logarithmic factor at most. On the
other hand, we show that the set $\{0,1\}^n\seq\Z^n$ possesses a dissociated
subset of size $\Ome(n\log n)$; since the standard basis of $\Z^n$ is a
maximal dissociated subset of $\{0,1\}^n$ of size $n$, the result just
mentioned is essentially sharp.
\end{abstract}

\maketitle

Recall, that \emph{subset sums} of a subset $\Lam$ of an abelian group are
group elements of the form $\sum_{b\in B} b$, where $B\seq\Lam$; thus, a
finite set $\Lam$ has at most $2^{|\Lam|}$ distinct subset sums.

A famous open conjecture of Erd\H os, first stated about 80 years ago (see
\refb{b} for a relatively recent related result and brief survey), is that if
all subset sums of an integer set $\Lam\seq[1,n]$ are pairwise distinct, then
$|\Lam|\le\log_2 n+O(1)$ as $n\to\infty$; here $\log_2$ denotes the base-$2$
logarithm. Similarly, one can investigate the largest possible size of
subsets of other ``natural'' sets in abelian groups, possessing the property
in question; say,

\begin{quote}
\emph{What is the largest possible size of a set $\Lam\seq\{0,1\}^n\seq\Z^n$
with all subset sums pairwise distinct?}
\end{quote}

In modern terms, a subset of an abelian group, all of whose subset sums are
pairwise distinct, is called \emph{dissociated}. Such sets proved to be
extremely useful due to the fact that if $\Lam$ is a maximal dissociated
subset of a given set $A$, then every element of $A$ is representable
(generally speaking, in a non-unique way) as a linear combination of the
elements of $\Lam$ with the coefficients in $\{-1,0,1\}$. Hence, maximal
dissociated subsets of a given set can be considered as its ``linear bases
over the set $\{-1,0,1\}$''. This interpretation naturally makes one wonder
whether, and to what extent, the size of a maximal dissociated subset of a
given set is determined by this set. That is,
\begin{quote}
\emph{Is it true that all maximal dissociated subsets of a given finite set
in an abelian group are of about the same size?}
\end{quote}

In this note we answer the two above-stated questions as follows.
\begin{theorem}\label{t:example}
For a positive integer $n$, the set $\{0,1\}^n$ (consisting of those vectors
in $\Z^n$ with all coordinates being equal to $0$ or $1$) possesses a
dissociated subset of size $(1+o(1))\,n\log_2n/\log_29$ (as $n\to\infty$).
\end{theorem}

\begin{theorem}\label{t:bound}
If $\Lam$ and $M$ are maximal dissociated subsets of a finite subset
$A\nsubseteq\{0\}$ of an abelian group, then
  $$ \frac{|M|}{\log_2(2|M|+1)} \le |\Lam|
        < |M| \, \big( \log_2(2M) + \log_2\log_2(2|M|) + 2 \big). $$
\end{theorem}

We remark that if a subset $A$ of an abelian group satisfies $A\seq\{0\}$,
then $A$ has just one dissociated subset; namely, the empty set.

Since the set of all $n$-dimensional vectors with exactly one coordinate
equal to $1$ and the other $n-1$ coordinates equal to $0$ is a maximal
dissociated subset of the set $\{0,1\}^n$, comparing
Theorems~\reft{example} and \reft{bound} we conclude that the latter is
sharp in the sense that the logarithmic factors cannot be dropped or
replaced with a slower growing function, and the former is sharp in the
sense that $n\log n$ is the true order of magnitude of the size of the
largest dissociated subset of the set $\{0,1\}^n$. At the same time, the
bound of Theorem~\reft{bound} is easy to improve given that the
underlying group has bounded exponent.

\begin{theorem}\label{t:exponent}
Let $A$ be finite subset of an abelian group $G$ of exponent $e:=\exp(G)$. If
$r$ denotes the rank of the subgroup $\<A\>$, generated by $A$, then for any
maximal dissociated subset $\Lam\seq A$ we have
  $$ r \le |\Lam| \le r\log_2 e. $$
\end{theorem}

We now turn to the proofs.
\begin{proof}[Proof of Theorem~\reft{example}]

We will show that if $n>(2\log_23+o(1))m/\log_2m$, with a suitable choice of
the implicit function, then the set $\{0,1\}^n$ possesses an $m$-element
dissociated subset. For this we prove that there exists a set
$D\seq\{0,1\}^m$ with $|D|=n$ such that for every non-zero vector
 $s\in S:=\{-1,0,1\}^m$ there is an element of $D$, not orthogonal to
$s$. Once this is done, we consider the $n\times m$ matrix whose rows are the
elements of $D$; the columns of this matrix form then an $m$-element
dissociated subset of $\{0,1\}^n$, as required.

We construct $D$ by choosing at random and independently of each other $n$
vectors from the set $\{0,1\}^m$, with equal probability for each vector to
be chosen. We will show that for every fixed non-zero vector $s\in S$, the
probability that all vectors from $D$ are orthogonal to $s$ is very small,
and indeed, the sum of these probabilities over all $s\in S\stm\{0\}$ is less
than $1$. By the union bound, this implies that with positive probability,
every vector $s\in S\stm\{0\}$ is not orthogonal to some vector from $D$.

We say that a vector from $S$ is of type $(m^+,m^-)$ if it has $m^+$
coordinates equal to $+1$, and $m^-$ coordinates equal to $-1$ (so that
$m-m^+-m^-$ of its coordinates are equal to $0$). Suppose that $s$ is a
non-zero vector from $S$ of type $(m^+,m^-)$. Clearly, a vector
$d\in\{0,1\}^m$ is orthogonal to $s$ if and only if there exists $j\ge 0$
such that $d$ has exactly $j$ non-zero coordinates in the
$(+1)$-locations of $s$, and exactly $j$ non-zero coordinates in the
$(-1)$-locations of $s$. Hence, the probability for a randomly chosen
$d\in\{0,1\}^m$ to be orthogonal to $s$ is
  $$ \frac{1}{2^{m^++m^-}}
           \sum_{j=0}^{\min\{m^+,m^-\}} \binom{m^+}j \binom{m^-}j
                 = \frac{1}{2^{m^++m^-}} \binom{m^++m^-}{m^+}
                       < \frac1{\sqrt{1.5(m^++m^-)}} \, . $$
It follows that the probability for \emph{all} elements of our randomly
chosen set $D$ to be simultaneously orthogonal to $s$ is smaller than
$(1.5(m^++m^-))^{-n/2}$.

Since the number of elements of $S$ of a given type $(m^+,m^-)$ is
$\binom{m}{m^++m^-}\binom{m^++m^-}{m^+}$, to conclude the proof it suffices
to estimate the sum
  $$ \sum_{1\le m^++m^-\le m}
       \binom{m}{m^++m^-}\binom{m^++m^-}{m^+} (1.5(m^++m^-))^{-n/2} $$
showing that its value does not exceed $1$.

To this end we rewrite this sum as
  $$ \sum_{t=1}^m \binom mt\, (1.5t)^{-n/2} \sum_{m^+=0}^t \binom t{m^+}
       = \sum_{t=1}^m \binom mt\, 2^t\, (1.5t)^{-n/2} $$
and split it into two parts, according to whether $t<T$ or $t\ge T$, where
$T:=m/(\log_2m)^2$. Let $\Sig_1$ denote the first part and $\Sig_2$ the
second part. Assuming that $m$ is large enough and
  $$ n>2\log_23\,\frac{m}{\log_2m}\,(1+\phi(m)) $$
with a function $\phi$ sufficiently slowly decaying to $0$, we have
  $$ \Sig_1 \le \binom mT 2^T 1.5^{-n/2}
                       < \lpr\frac{9m}T\rpr^T 1.5^{-n/2}
                                         = (3\log_2m)^{2T} 1.5^{-n/2} , $$
whence
 $$ \log_2\Sig_1 < \frac{2m}{(\log_2m)^2}\,\log_2(3\log_2m) -
                 \log_23\log_2 1.5\,\frac{m}{\log_2m} \, (1+\phi(m)) < -1, $$
and therefore $\Sig_1<1/2$. Furthermore,
  $$ \Sig_2 \le T^{-n/2} \sum_{t=1}^m \binom mt 2^t < T^{-n/2} 3^m, $$
implying
\begin{align*}
  \log_2\Sig_2
    &< m\log_23 - (\log_2m - 2\log_2\log_2 m)\, \log_23
                                          \, \frac{m}{\log_2m}\,(1+\phi(m)) \\
    &= m\log_23 \lpr \frac{2\log_2\log_2m}{\log_2m}\,(1+\phi(m))-\phi(m) \rpr \\
    &< -1.
\end{align*}
Thus, $\Sig_2<1/2$; along with the estimate $\Sig_1<1/2$ obtained above, this
completes the proof.
\end{proof}

\begin{proof}[Proof of Theorem~\reft{bound}]
Suppose that $\Lam,M\seq A$ are maximal dissociated subsets of $A$. By
maximality of $\Lam$, every element of $A$, and consequently every element of
$M$, is a linear combination of the elements of $\Lam$ with the coefficients
in $\{-1,0,1\}$. Hence, every subset sum of $M$ is a linear combination of
the elements of $\Lam$ with the coefficients in $\{-|M|,-|M|+1\longc |M|\}$.
Since there are $2^{|M|}$ subset sums of $M$, all distinct from each other,
and $(2|M|+1)^{|\Lam|}$ linear combinations of the elements of $\Lam$ with
the coefficients in $\{-|M|,-|M|+1\longc |M|\}$, we have
  $$ 2^{|M|} \le (2|M|+1)^{|\Lam|}, $$
and the lower bound follows.

Notice, that by symmetry we have
  $$ 2^{|\Lam|} \le (2|\Lam|+1)^{|M|}, $$
whence
  $$ |\Lam| \le |M| \log_2(2|\Lam|+1). \eqno{(\ast)} $$

Observing that the upper bound is immediate if $M$ is a singleton (in which
case $A\seq\{-g,0,g\}$, where $g$ is the element of $M$, and therefore every
maximal dissociated subset of $A$ is a singleton, too), we assume $|M|\ge 2$
below.

Since every element of $\Lam$ is a linear combination of the elements of $M$
with the coefficients in $\{-1,0,1\}$, and since $\Lam$ contains neither $0$,
nor two elements adding up to $0$, we have $|\Lam|\le(3^{|M|}-1)/2$.
Consequently, $2|\Lam|+1\le 3^{|M|}$, and using ($\ast$) we get
  $$ |\Lam| \le |M|^2\log_23. $$
Hence,
  $$ 2|\Lam|+1 < |M|^2\log_29 + 1 < 4|M|^2, $$
and substituting this back into ($\ast$) we obtain
  $$ |\Lam| < 2|M|\log_2(2|M|) . $$
As a next iteration, we conclude that
  $$ 2|\Lam|+1 < 5|M|\log_2(2|M|), $$
and therefore, by ($\ast$),
  $$ |\Lam| \le |M| \big( \log_2(2|M|)
                                  + \log_2\log_2(2|M|) + \log_2(5/2) \big). $$
\end{proof}

\begin{proof}[Proof of Theorem~\reft{exponent}]
The lower bound follows from the fact that $\Lam$ generates $\<A\>$, the
upper bound from the fact that all $2^{|\Lam|}$ pairwise distinct subset sums
of $\Lam$ are contained in $\<A\>$, whereas $|\<A\>|\le e^r$.
\end{proof}

\bigskip
We close our note with an open problem.
\begin{quote}
\emph{For a positive integer $n$, let $L_n$ denote the largest size of a
dissociated subset of the set $\{0,1\}^n\seq\Z^n$. What are the limits
  $$ \liminf_{n\to\infty} \frac{L_n}{n\log_2n}\ \text{and}
       \ \limsup_{n\to\infty} \frac{L_n}{n\log_2n} \,? $$}
\end{quote}
Notice, that by Theorems~\reft{example} and \reft{bound} we have
  $$ 1/\log_29 \le \liminf_{n\to\infty} \frac{L_n}{n\log_2 n}
               \le \limsup_{n\to\infty} \frac{L_n}{n\log_2n} \le 1. $$

\bigskip

\end{document}